\documentclass{article}

\usepackage{arxiv}

\usepackage[utf8]{inputenc} 
\usepackage[T1]{fontenc}    
\usepackage{hyperref}       
\usepackage{url}            
\usepackage{booktabs}       
\usepackage{amsfonts}       
\usepackage{nicefrac}       
\usepackage{microtype}      
\usepackage{lipsum}		
\usepackage{graphicx}
\usepackage{natbib}
\usepackage{doi}

\title{Polynomial algorithms for Simultaneous Unitary Similarity and Equivalence}

\author{ \href{https://orcid.org/0000-0000-0000-0000}{\includegraphics[scale=0.06]{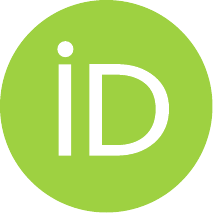}\hspace{1mm}Harikrishna VJ} \\
	Mysuru\\
	India\\
	\texttt{harikrishna.vj@gmail.com} \\
	\And
	\href{https://orcid.org/0000-0000-0000-0000}{\includegraphics[scale=0.06]{orcid.pdf}\hspace{1mm}Vittal Rao} \\
	Bengaluru\\
	India\\
	\texttt{vittalrao.14@gmail.com} \\
		\And
	\href{https://orcid.org/0000-0000-0000-0000}{\includegraphics[scale=0.06]{orcid.pdf}\hspace{1mm}K.R Ramakrishnan} \\
	Bengaluru\\
	India\\
	\texttt{krr2504@gmail.com} \\
}

\renewcommand{\shorttitle}{\textit{arXiv} Template}

\hypersetup{
pdftitle={A template for the arxiv style},
pdfsubject={q-bio.NC, q-bio.QM},
pdfauthor={David S.~Hippocampus, Elias D.~Striatum},
pdfkeywords={First keyword, Second keyword, More},
}

\begin{document}
\maketitle

\begin{abstract}
	We present a polynomial algorithm to solve the Simultaneous Unitary Similarity (S.U.S) problem i.e to find a Unitary $U$ if it exists s.t $UA_lU^*=\mathcal{B}_l$, $l \in \{1,\cdots,p\}$, where $(\mathcal{A}_l,\mathcal{B}_l)$'s are $n$x$n$ matrix $2$-tuples over $C$ the field of Complex numbers. We show that the algorithm naturally extends to solve the Simultaneous Unitary Equivalence (S.U.Eq) problem which is to find if there exist Unitary $U,V$ s.t $UA_lV^*=\mathcal{B}_l$, where $(\mathcal{A}_l,\mathcal{B}_l)$'s are Complex Rectangular $m$x$n$ matrix $2$-tuples.
	
	To solve the S.U.S problem, we observe that when $U$ is a diagonal matrix i.e $U=diag(u_1,\cdots,u_n)$ we can solve for the scalars $u_1,\cdots,u_n$ by connecting them using `paths' in a graph whose edges are non-zero entries of $(\mathcal{A}_l,\mathcal{B}_l)$'s . Inspired by this, in the general case when $U$ is block-diagonal i.e $U=diag(U_1,\cdots,U_d)$ we define the `Solution-form' for the collection $\{(\mathcal{A}_l,\mathcal{B}_l)\}_{l=1}^p$ and partition the $d$ unitary blocks using paths of unitary sub-matrices to solve for $U$. When the collection is not in Solution-form we show that by diagonalizing a Hermitian or a Normal matrix related to the sub-matrices an equivalent problem can be set-up where $U$ has more blocks than $d$. Thus in every iteration the problem either simplifies or gets solved altogether.
	
	The solution extends naturally to solve the S.U.Eq problem where both $U$ and
	$V$ are visualized as block-diagonal and in each step either one of them `blocks further' or we find paths in a bi-graph of vertices defined by blocks of $U,V$ to partition the vertices and solve the problem. \\
	The algorithms have a polynomial run time complexity of $O(pn^4)$.\\
	The problem has application in Quantum Computation, Quantum Gate Design \& Quantum Simulation.
	
	For a collection $\{(\mathcal{A}_l)\}_{l=1}^p$ of $p$ matrices if the steps of the algorithm are followed until we get to the Solution form then the salient results of each step act as invariant factors or Canonical features thus classifying the collection upto Simultaneous Unitary Similarity.
\end{abstract}


\section{Introduction}
The analysis of Similarity of Matrices over Fields, as well as Integral Domains which are not Fields, is a classical problem in Linear Algebra and has
received considerable attention. A related problem is that of Simultaneous
Similarity of Matrices. The many interesting algebraic questions and the history of the problems is discussed by Shmuel Friedland \citep{frdlnd1983SimultSimil}. 

A special case of this
problem is that of Simultaneous Unitary Similarity (S.U.S) of matrices which is 
to find a Unitary $U$ if it exists s.t $UA_lU^*=\mathcal{B}_l$, $l \in \{1,\cdots,p\}$, where $(\mathcal{A}_l,\mathcal{B}_l)$'s are given $n$x$n$ matrix $2$-tuples over $C$ the field of Complex numbers.
Another related problem is referred to as
the Simultaneous Unitary Equivalence (S.U.Eq) problem where we find Unitary $U,V$ if they exist s.t $UA_lV^*=\mathcal{B}_l$, where $(\mathcal{A}_l,\mathcal{B}_l)$'s are $p$ given $2$-tuples of $m$x$n$ matrices over $C$.

The problem has application in Quantum Computation, Quantum Gate Design and related areas since any state in Quantum Mechanics is essentially a semi-definite matrix and
evolution is Unitary Similarity Operation. So whereever questions of Simultaneous Evolution of certain source Quantum states to given destination states come up the 
S.U.S problem finds an application.\citep{Jing2014SimultUnitaryApplctn,Chefles2000SimultUnitry,Chefles2004SimultUnitry}

The problem can be treated as the problem of classifying matrices upto Simultaneous Unitary/Orthogonal Similarity as in Friedland's \citep{frdlnd1983SimultSimil}\cite{Shapiro1991Survey} work where the characteristic polynomials of $\Sigma^p_{j=1} \mathcal{A}_jx^j$ and $\Sigma^p_{j=1} \mathcal{B}_jx^j$ are used to describe the invariant factors. Radjavi \citep{hydrrdjvi1968SimultSimil} describes an algorithm in which invariant factors are found based on diagonalization of a sub-matrix in each step and by adding matrices to the collection of $2$-tuples whose Simultaneous Unitary Similarity is being checked. The algorithm leads to a Canonical form with a complexity which is $3^n$ times the cost of each step.While these approaches solve the problem theoretically comprehensively they do not lead to an efficient algorithm to solve the S.U.S and S.U.Eq problems.

Recently to solve the S.U.S problem Gerasimova et al \citep{Gerasimova2013SimultEq} construct words of products of $\mathcal{A}_l$'s and $\mathcal{B}_l$'s whose traces must agree for them to be in the same equivalence class. This is along the lines of Specht's theorem \citep{Specht1940Trace} for the non-simultaneous version of this problem which is to find if there exists a $U$ s.t $U\mathcal{A}U^*=\mathcal{B}$ over Complex numbers. Although it is proved that the length of the words is bounded \cite{Pappacena1997UpprBnd} it still does not lead to an efficient implementation as the number of words would be prohibitively large. Jing \citep{Jing2015SimultEq} describes a similar trace matching algorithm for the S.U.Eq problem. 
From the point of view of computational complexity generally speaking any problem involving Classification of matrices upto similarity is considered to be a wild problem. \citep{Segei1998wild}
 
We present algorithms to solve the S.U.S and S.U.Eq problems which have polynomial running time complexity in terms of size of matrices and the number of $2$-tuples.
Our algorithm is based on the observation that the problem is simplest when 
the unitary $U$ we are looking for is diagonal i.e $U=diag(u_1,\cdots,u_n)$. 
In this case each modulo one scalar can be uniquely written in-terms of the other by using non-zero elements in $\mathcal{A}_l$'s and $\mathcal{B}_l$'s. The problem has the structure of a graph with edges being the non-zero elements and `paths' of these non-zero elements partitioning the $n$ vertices. To generalize this we define a Solution-form for the collection of $2$-tuples which resembles the simplest case when the unitary $U$ is $diag(U_1,U_2,\cdots,U_d)$. The `paths' that partition $\{1,2,\cdots,d\}$ and inter-relate the Unitary blocks are formed using the sub-matrices of $\mathcal{A}_l$'s and $\mathcal{B}_l$'s which are non-zero multiples of Unitary. If the matrix is not in Solution-form we show that a Hermitian matrix or a Normal matrix formed using the sub-matrices can be diagonalized to set-up an equivalent problem where $U$ has more blocks than $d$ thus moving a step closer to the simplest case and solving the problem in a maximum of $n$ steps.     
The algorithm naturally extends to solve the S.U.Eq problem where in each step we work with Unitary $U=diag(U_1,U_2,\cdots,U_d)$ and $V=diag(V_1,V_2,\cdots,V_f)$. In the Solution-form we form a bi-graph with vertices $\{1_R,\cdots,d_R\}$ $\cup$ $\{1_C,\cdots,f_C\}$ ($R,C$ indicate rows and columns respectively) and partition it to identify the $U$ and $V$ that should solve the problem. If we do not reach Solution-form the problem is simplified with either $U$ or $V$ having more blocks than before. Thus the problem is solved in a maximum of $m+n$ steps.

The algorithms are discussed in detail in section \ref{sec:Algorithms}.
The complexities  of checking for Solution form and setting up the Equivalent problem dominate over other costs leading to a polynomial runtime cost of $O(pn^4)$ as discussed in \ref{sec:comp_compl}
Although the diagonalization step is numerically stable the accuracy depends on the eigen-gap in the matrix being diagonalized. Our implementation, related stability and accuracy issues are discussed in section \ref{subsec:numerical}.
We also discuss in section \ref{subsec:numerical} the `Canonical features' that
could be retained in each step if we work with a collection of matrices $\{A_l\}_{l=1}^p$ instead of collection of $2$-tuples

It is evident that the research in this field has a rich history and most of the work published on this uses techniques from advanced Linear Algebra, Algebraic Geometry and Representation Theory \citep{frdlnd1983SimultSimil,Gerasimova2013SimultEq}.
Compared to these our approach is closer to the Radjavi's \citep{hydrrdjvi1968SimultSimil} approach to solve the $S.U.S$ problem.
 An important difference is that we use graph based connectivity to simplify the problem in
each iteration without adding to the inital collection of matrices. To get the solution we work with the `paths' in the graphs whose number is $\le n$ leading to an algorithm of polynomial time complexity. We extend the technique to solve the $S.U.Eq$ problem as well. 
Another difference being we work with the original collection it self without having to get Hermitian pairs for each matrix in the collection. 
Compared to methods that use techniques from algebraic geometry and Representation theory our work is easy to grasp as it uses constructs such as diagonalizing a Hermitian or a Normal matrix, partitioning a graph which even a keen beginner in the field can understand and implement.
Linking of a graph, edges and connectedness to solve problems involving matrices is not new and
has been used before mainly to check irreducibility of matrices with applications in spectral theory and for checking if a matrix can be triangulized\citep{meyer2000irdcblty}. In the context of solving the S.U.S and S.U.Eq problems our approach of using an undirected graph, a bi-graph with paths helping us to identify the form needed to solve the problem is novel.

\section{Algorithms for S.U.S and S.U.Eq problems}
\label{sec:Algorithms}
In this section we describe algorithms to solve the S.U.S problem
and then show that the algorithm can be naturally extended to solve
the S.U.Eq problem.  

The S.U.S problem can be stated as:
Given a collection of $p$ $n$x$n$ Complex $2$-tuples $\{(\mathcal{A}_l,\mathcal{B}_l)\}_{l=1}^{p}$ 
find a Unitary $U$ if it exists s.t $U\mathcal{A}_lU^*=\mathcal{B}_l$.

The S.U.Eq problem:
Given $p$, $m$x$n$ Complex $2$-tuples $\{(\mathcal{A}_l,\mathcal{B}_l)\}_{l=1}^{p}$, solve for Unitary $U,V$ s.t $U\mathcal{A}_lV^*=\mathcal{B}_l$ 

\subsection{The simplest case, diagonal $U$ and paths of non-zero elements}
\begin{figure}
	\centering
	\includegraphics[width=100mm]{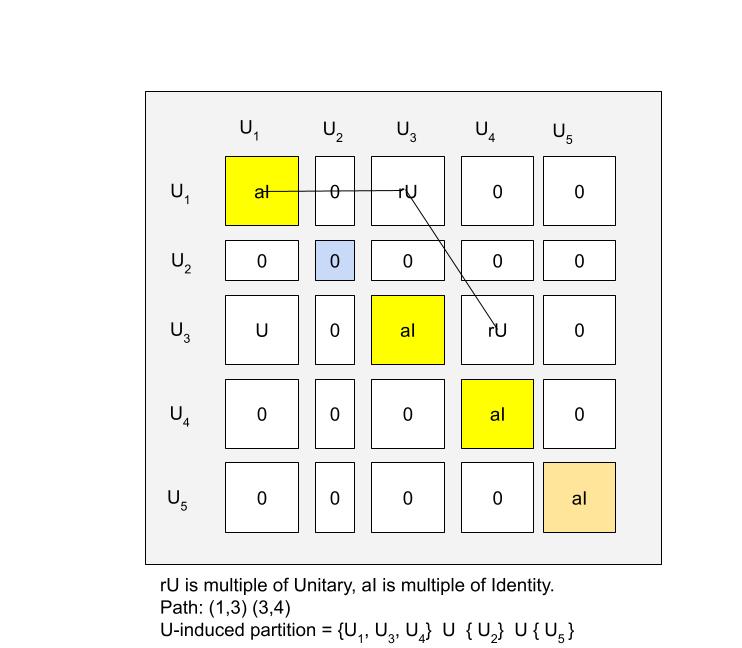}
	\caption{Pictorial depiction of S.U.S's Solution form with paths and partition of blocks in $U$ }
	\label{fig:fig1}
\end{figure}

The algorithm is constructed based on the ideas that help to solve the problem
when the $U$, $U=(u_1,u_2,\cdots,u_n)$ is diagonal with $u_{i=1}^n$
being modulo $1$ complex numbers to be solved for. Let $n=4$, and $p=1$
i.e a single $2$-tuple $(A,B)$ with $[A]_{ij}=a_{ij}$ and $[B]_{ij}=b_{ij}$. 
We are solving for modulo $1$ $u_i,u_j$ in $u_ia_{ij}u_j^*$=$b_{ij}$. 
$i,j \in \{1,2,3,4\}$. W.l.o.g let $u_1=1$ assume $a_{12}$ is non-zero
then using $u_1a_{12}u_2^*$=$b_{12}$ we can solve for $u_2$. Now the question
arises, which are all the $u_i$'s that can be solved if $u_1$ is known. Naturally the $u_i$'s that are 'connected' to $u_1$ by non-zero elements 
can be solved. For eg: Suppose $a_{13}$ and $a_{34}$ are non-zero
then we can consider a product of corresponding equations to give :
\begin{eqnarray*}
	u_1a_{13}u_3^*u_3a_{34}u_4^*=b_{13}b_{34}\\ \\
	u_1a_{13}a_{34}u_4^*=b_{13}b_{34}
\end{eqnarray*} 
Hence $u_4$ gets related to $u_1$ through the `paths' of non-zero elements along the way and can be written in terms of $u_1$. In general all $u_i$'s
that are related to $u_1$ via these paths are related to $u_1$ and get solved in terms of $u_1$. The $u_i$'s that are not related $u_1$ form another
`class' and they can be solved in terms of of a representative of their class
using the paths. Thus the graph with vertices $\{1,\cdots,n\}$ and the notion of paths defining relation between vertices (Figure \ref{fig:fig1}) which gives rise to equivalence classes of `related' variables is the structure that gives the solution. 
Inspired by this idea we think of the general case as solving for $U=diag(U_1,\cdots,U_d)$ where $U$ is block-diagonal ($U=diag(U_1)$ in the first iteration with $U_1=U$) and consider the submatrices in the collection of 
$2$-tuples formed by partitioning the matrices based on sizes of $U_1,\cdots,U_d$.
This allows us to identify a form of the collection of $2$-tuples when a solution can be obtained as in the simplest case using the paths and equivalence classes of $U_i$'s. We refer to this as the Solution form. 
The non-zero elements of the simplest case which define the 'paths' get generalized to non-zero multiples of unitary and $U^{sol}$ that should solve the problem
is formed using products along `paths' of these non-zero elements.
We show that when the collection is not in Solution-form we can identify 
a sub-matrix in the collection, diagonalize and find an `Equivalent' problem
with $U$ having more blocks than before. Thus in each step
we either move a step closer to the simplest case or solve the problem.
 
The same ideas naturally extend to solve the S.U.Eq problem:
where we solve for $U=diag(U_1,\cdots,U_d)$ and $V=diag(V_1,\cdots,V_f)$ using the the bi-graph of vertices $\{1_R,\cdots,d_R\}$  $\cup$ $\{1_C,\cdots,f_C\}$ ( $R$,$C$ indicate row and column respectively ) to define paths through the sub-matrices of the collection and get the Solution-form.

In the subsections to follow we describe the flow of the algorithm, define the Solution form, the associated graph, paths, the eventual solution, describe how an Equivalent problem is constructed and prove that the algorithms work.

We first explain the Algorithm for the S.U.S problem in detail and then show 
that the results with appropriate changes can be used to solve the S.U.Eq problem as well.

\subsection {Algorithm flow}
\label{subsec:flow}
The algorithm runs for a maximum of $n$ iterations and following is the computation involved in each iteration.

\begin{tabbing}
	\textbf{Input:} \= Collection $\{(\mathcal{A}_l,\mathcal{B}_l)\}_{l=1}^p$, \\
					\>Solve for Unitary  $U=diag(U_1,\cdots,U_d)$ s.t $U\mathcal{A}_lU^*=\mathcal{B}_l$ \\ \\

	\textbf{If:} \> Collection in Solution-form :\\
		\>Using the $U$-induced graph of $\{(\mathcal{A}_l,\mathcal{B}_l)\}_{l=1}^p$ and the paths that partition $\{1,\cdots,d\}$\\ 
		\>compute $U^{sol}=diag(U^{sol}_1,\cdots,U^{sol}_d)$\\ \\
		
	\textbf{Else:}
		\>Diagonalize and Compute an equivalent problem with  \\
		\>$U$ having the structure $U=diag(U_1,\cdots,U_{d^{'}})$, $d^{'}>d$	
\end{tabbing}

In each iteration of the algorithm we check if the collection is in Solution form and if not in Solution form we construct an Equivalent problem that gets us a step closer to the simplest case solving the problem in a maximum of $n$ steps. In the ensuing subsections we describe the terms required to understand the
Solution form and the construction of the Equivalent problem. The main results which give us the $U^{sol}$ in solution form and the Equivalent problem are stated and proved as Theorem 1 and Theorem 2.

\subsection{Solution form}
To get to $U^{sol}$ we partition the matrices of the collection along rows and columns based on the sizes of the blocks in $U$ and define a $U$-induced graph with edges defined by non-zero submatrices. The products of these sub-matrices along the `paths' give $U^{sol}$. The main result is stated and proved in Theorem 1. 
\subsubsection {U-Induced Partition \& Solution form }
The partition of the given collection based on structure of $U$ is defined as 
follows:

\textbf{Definition 1:} $U$-Induced Partition of the given collection of matrices $\{(\mathcal{A}_l,\mathcal{B}_l)\}_{l=1}^{p}$ where $U=diag(U_1,\cdots,U_d)$.
Suppose the size of the $d$ diagonal blocks of $U$ is $n_1,n_2,\cdots,n_d$
then the partition of the matrices of the collection into first $n_1$ rows and then the next $n_2$ rows upto the last $n_d$ rows ($\Sigma_{l=1}^p n_i=n$) and similarly along columns is the $U$-Induced Partition of the given collection. 

The Pre-Solution form is defined first for the given collection based on the sub-matrices in the $U$-induced partition of the given collection and once the associated graph is found a the structure that gives the Solution form is defined.

\textbf{Definition 2:} Pre-Solution form. Given that the unitary $U$ being solved for is $U=diag(U_1,\cdots,U_d)$ the collection of $2$-tuples $\{(\mathcal{A}_l,\mathcal{B}_l)\}^p_{l=1}$ is said to be in Pre-Solution form if the $U$-Induced Partition of the $\{(\mathcal{A}_l,\mathcal{B}_l)\}^p_{l=1}$ has the following structure:
\begin{enumerate}
	\item All diagonal sub-matrices are multiples of Identity 
	\item All non-square sub-matrices are zero.
	\item All off diagonal square matrices are multiples of Unitary.
\end{enumerate}
The following Lemma relates the sub-matrices of the $A_l$ with those of $B_l$ in the Pre-Solution form. 

\textbf{Lemma 1:} Given $U=diag(U_1,\cdots,U_d)$ and a collection $\{(\mathcal{A}_l,\mathcal{B}_l)\}^p_{l=1}$ in Pre-Solution form. If the S.U.S problem has a solution then the following are true:
\begin{enumerate}
	\item All diagonal sub-matrices are multiples of Identity with $A_{ii}=\alpha_{ii}I=B_{ii}$. $i,j \in \{1,\cdots,d\}$.
	\item Non-square submatrices are zeros i.e $A_{ij}=0=B_{ij}$.
	\item All off-diagonal square matrices are multiples of Unitary, i.e $A_{ij}A^*_{ij}=a_{ij}I=B_{ij}B^*_{ij}$, $a_{ij}\ge 0$
\end{enumerate}

\textbf{Proof:}
\begin{enumerate}
	\item In Solution form given there exists a solution, if $A_{ii}=\alpha_{ii}I_{n_i}$, $U_i$ is $n_i$x$n_i$.  $B_{ii}=U_iA_{ii}U^*_i=U_i\alpha_{ii}I_{n_i}U^*_i=\alpha_{ii}I_{n_i}$.
	\item If $A_{ij}=0$ then $B_{ij}=U_iA_{ij}U^*_j$ is also $0$.
	\item $B_{ij}B^*_{ij}$=$U_iA_{ij}U^*_jU_jA^*_{ij}U^*_i=U_iA_{ij}A^*_{ij}U^*_i=U_iaI_{n_i}U^*_i=aI_{n_i}$, Since $A_{ij}A^*_{ij}$ is positive semi-definite $a\ge 0$.
\end{enumerate}

If any of these conditions is not met either the collection is not in Pre-Solution form
or it is evidence that the $A_l$'s are not unitarily similar to $B_l$'s.

\subsubsection{ $U$-induced graph and paths that partition $U$ }
A collection in Pre-Solution form allows us to define a graph, its edges and paths and when Solution form is defined it resembles the simplest 
case (diagonal $U$) where we solve for modulus-$1$ scalars by connecting the variables using paths of non-zero entries. Solution form allows us to solve for the blocks in $U$ by using the paths of non-zero elements which are multiples of Unitary. To prove the above we need the following definitions.

\textbf{Definition 3:}
Given $U=diag(U_1,\cdots,U_d)$, and a collection $\{(\mathcal{A}_l,\mathcal{B}_l)\}_{l=1}^{p}$ in Pre-Solution form the $U$-Induced Graph of the given collection of matrices is a graph with vertices $i,j \in $ $\{1,\cdots,d\}$, where an edge between vertices $i,j$ is said to exist if  $A_{ij}$ ($B_{ij}$) or $A_{ji}$ ($B_{ji}$) is non-zero where $A_{ij}$, $B_{ij}$ \& $A_{ji}$ , $B_{ji}$ are submatrices in the $U$-Induced Partition of the given collection.
To simplify the notation of the submatrices  we do not include as part of notation the matrix number $l \in$ $\{1,\cdots,p\}$ to which the non-zero sub-matrices belong.

\textbf{Lemma 2:} Given a $U$-induced graph of a collection $\{(\mathcal{A}_l,\mathcal{B}_l)\}_{l=1}^{p}$ where $U$ is Unitary with the form $U=diag(U_1,\cdots,U_d)$  the path that connects any two vertices of this graph defines an equivalence relation on the set of vertices i.e $\{1,\cdots,d\}$ and hence partitions the set of vertices. 

\textbf{Proof:} The fact that $U=diag(U_1,\cdots,U_d)$ is Unitary can be written
as $U\mathcal{D}U^*=\mathcal{D}$ with $\mathcal{D}=diag(\lambda_1I_{n_1},\cdots,\lambda_dI_{n_d})$ , where the $\lambda$'s
can all be assumed to be greater than zero. This is just to capture the fact that
$U$ is unitary and has the given block diagonal structure.
The $2$-tuple $(\mathcal{D},\mathcal{D})$ can be treated as being part of the given collection. Since $\lambda_i$, $i \in \{1,\cdots,d\}$ is non-zero each $i$ is connected to itself by
the path $(i,i)$ with $\lambda_i I_{n_i}$ being the non-zero diagonal element in $\mathcal{D}$. Hence the relation defined by the existence of a path between vertices
in the $U$-Induced Graph is reflexive.\\
To prove Symmetry suppose we have a path from $i$ to $j$, this path is formed by sequence of edges in the graph corresponding to which there are non-zero sub-matrices in the collection. If we reverse this sequence then we have a path from $j$ to $i$ since the graph is not a directed graph.\\
To prove Transitivity, Suppose there is a path from $i$ to $j$ and another from $j$
to $k$ the concatenation of these paths will serve as a path from $i$ to $k$ hence
the relation is Transitive. In Summary the relation is an Equivalence relation on the set of vertices and hence partitions the set of vertices i.e $\{1,\cdots,d\}$.

The following definition of products along the paths in the partition of the graph helps us define Solution form.

\textbf{Definition 4:} $A^{pth}_{c(i)i}$, $B^{pth}_{c(i)i}$ \&  $pr(A_{ij})$, $pr(B_{ij})$ .  The class to which index $i$ belongs is denoted by $c(i)$.
Given the $U$-Induced Partition of a collection of $2$-tuples in Pre-Solution form,
$A^{pth}_{c(i)i}$ and $B^{pth}_{c(i)i}$ denote product of submatrices
along the path from $c(i)$ to $i$, formed as explained below.
Suppose the path from $c(i)$ to $i$ in $A_l$'s is represented by the sequence of vertices $v_1,v_2,v_3,v_4$ with $v_1=c(i)$ and $v_4=i$ and the sequence of edges is $(i_1,j_1),(i_2,j_2),(i_3,j_3)$. 
\begin{eqnarray*}
	A^{pth}_{c(i),i}=A^{s(0)}_{i_0j_0}A^{s(1)}_{i_1j_1}A^{s(2)}_{i_2j_2}\\
	B^{pth}_{c(i),i}=B^{s(0)}_{i_0j_0}B^{s(1)}_{i_1j_1}B^{s(2)}_{i_2j_2}
\end{eqnarray*}
Here $s$ is a function which takes value $1$ or $-1$ based on whether the 
$k^{th}$ edge is $(v_k,v_{k+1})$ or $(v_{k+1},v_{k})$ since either of these
could be non-zero to define an edge. Let us look at an example with $c(i)=1$, 
$i=4$ with the sequence of vertices along the path being $1,2,3,4$ and the edge sequence being $(2,1),(2,3),(4,3)$ then $A^{pth}_{14}$ is given by:
\begin{equation}
A^{pth}_{1,4}=A^{-1}_{21}A_{23}A^{-1}_{43}
\end{equation}
$s(1)$ is $-1$ since $(2,1)$=$(v_2,v_1)$ and not $(v_1,v_2)$. Similarly $s(3)$ is $-1$ since $(4,3)=(v_4,v_3)$, whereas $s(2)$ is
1 as $(2,3)$ is $(v_2,v_3)$. The inverses exist since
the submatrices involved are non-zero multiples of unitary by Lemma 3 (3rd part).

Given ${\{(\mathcal{A}_l,\mathcal{B}_l)\}}_{l=1}^{p}$ in Pre-Solution form using the paths in the $U$-Induced Partition, with $i,j$ in the same class of the partition, we define $pr(A_{ij})$, $pr(B_{ij})$ using $A^{pth}_{c(i)i}$ ($B^{pth}_{c(i)i}$) and $A^{pth}_{c(j)j}$ ($B^{pth}_{c(j)j}$), as follows:
\begin{eqnarray*}
	pr(A_{ij})&=&A^{pth}_{c(i)i} A_{ij} (A^{pth}_{c(i)j})^{-1} \\
	pr(B_{ij})&=&B^{pth}_{c(i)i} B_{ij} (B^{pth}_{c(i)j})^{-1}
\end{eqnarray*}  
Since $i,j$ belong to the same equivalence class $c(i)=c(j)$ and the above definition applies to $A_{ij}$ and $B_{ij}$ in all the matrices of the collection.

\textbf{Lemma 3:} Given $U$-Induced Partition of a collection $\{(\mathcal{A}_l,\mathcal{B}_l)\}_{l=1}^{p}$ in Pre-Solution form, if the collection is Simultaneously Unitarily Similar then:
\begin{eqnarray*} 
 U_{c(i)}pr(A_{ij})U^*_{c(i)}=pr(B_{ij}) \\
\end{eqnarray*}
\textbf{Proof:}
We prefer to prove by taking an example as it simplifies the notation.
W.l.o.g let us consider an example where $3,4$ belongs to the class $\{1,3,4\}$ in the $U$-induced partition. The path from $1$ to $3$ being $(1,3)$  and the path from $1,4$ being $(1,3)(4,3)$. Here $i=3,j=4$ and $c(i)=1=c(j)$.
\begin{eqnarray*}
	U_1pr(A_{34})U^*_1&=& U_1A^{pth}_{13} A_{34} (A^{pth}_{14})^{-1} U^*_1 \\
	&=&U_1 A_{13} A_{34} (A_{13}A^{-1}_{43})^{-1} U^*_1 \\
	&=&U_1 A_{13} U_3U^*_3 A_{34} U^*_4U_4 A_{43}U_3 U^*_3A^{-1}_{13} U^*_1\\
	&=&B_{13} B_{34} (B_{13}B^{-1}_{43})^{-1} \\
	&=& B^{pth}_{13} B_{34} (B^{pth}_{14})^{-1} \\
	&=& pr(B_{34})
\end{eqnarray*} 
Hence the proof.

We have set out to get to a stage that best resembles the simplest case described in Section 2.1, where modulo-$1$ scalars
to be solved for in the diagonal $U$ are written in terms of a variable which represents the class and is a free variable. Here $U_{c(i)}$ is the representative of a class and
hence it has to work like a free variable. The only way $U_{c(i)}$
can be a free variable and at the same time satisfy Lemma 3 is
if $pr(A_{ij})$ is a multiple of identity. We add this condition to
the definition of Pre-Solution form to get Solution form and
define $U^{sol}$ which solves the problem for a collection in Solution form.

\textbf{Definition 5:} Solution form \& $U^{sol}$. 
Given a collection $\{(\mathcal{A}_l,\mathcal{B}_l)\}_{l=1}^{p}$ in Pre-Solution form and given the $U$-Induced Graph and it's partition, the collection is said to be in Solution form if $pr(A_{ij})$, $pr(B_{ij})$ for $i,j \in \{1,\cdots,d \}$ are multiples of identity i.e :
\begin{eqnarray*}
 pr(A_{ij})=\beta_{ij}I=pr(B_{ij}) \\
\end{eqnarray*}
If $pr(A_{ij})=\beta_{ij}I$ then by Lemma 3, $pr(B_{ij})$ = $U_{c(i)}pr(A_{ij})U^*_{c(i)}=\beta_{ij}I$ 

In Solution form we can define the $U$ that solves the problem as $U^{sol}=(U^{sol}_1,\cdots,U^{sol}_d)$ using the paths in the $U$-Induced Graph of the collection. Let $c(i)$, $i \in \{1,\cdots,d \}$ denote the class representative of vertex $i$, using the path from $c(i)$ to $i$, $U^{sol}_i$
is defined in terms of $U_{c(i)}$ as:
\begin{eqnarray*}
	U^{sol}_i=(B^{pth}_{c(i)i})^{-1} U_{c(i)} A^{pth}_{c(i)i}
\end{eqnarray*}

Following theorem shows that if at all there is a solution to the S.U.S problem then $U^{sol}$ is a solution.

\textbf{Theorem 1:} Let the $U$ we are looking to solve in the S.U.S problem be $U=diag(U_1,\cdots,U_d)$. If the given collection ${\{(\mathcal{A}_l,\mathcal{B}_l)\}}_{l=1}^{p}$ is in Solution form then the
given S.U.S problem has a solution i.f.f $U^{sol}$ as given in Definition 4
is a solution.

\textbf{Proof:}
Given the collection and the $U$-Induced Partition of the collection as per Definition 3, we can get the $U$-Induced Graph of the collection and partition
the indices $\{1,\cdots,d\}$ into equivalence classes as proved in Lemma 2. Then as denoted in Definition 4, 
we can associate products with the
paths from the class representative of a vertex $i$ in the set of vertices $\{1,\cdots,d\}$ denoted by $c(i)$ to $i$. W.l.o.g let us consider an example
with $i=4$, $c(i)=1$ the sequence of vertices $1,2,3,4$ with the path being
$(1,2),(3,2),(3,4)$ i.e $s(1)=1$, $s(2)=-1$ and $s(3)=1$. The path-products for vertex $4$ would be:
\begin{eqnarray*}
	A^{pth}_{1,4}&=&A_{12}A^{-1}_{32}A_{34}\\
	B^{pth}_{1,4}&=&B_{12}B^{-1}_{32}B_{34}
\end{eqnarray*}
Consider $U_{1}A^{pth}_{1,4}U_{4}^*$, 
\begin{eqnarray*}
 	U_{1}A^{pth}_{1,4}U^*_4&=&U_1A_{12}A^{-1}_{32}A_{34}U^*_4\\
 	&=&U_1A_{12} U^*_2U_2 A^{-1}_{32} U^*_3U_3 A_{34}U^*_4 \\
 	&=&U_1A_{12}U^*_2 (U_3 A_{32} U^*_2)^{-1} U_3A_{34}U^*_4 \\
 	&=&U_1A_{12}U^*_2 (U_3 A_{32} U^*_2)^{-1} U_3A_{34}U^*_4 \\
 	&=&B_{12}B^{-1}_{32}B_{34}\\
 	&=&B^{pth}_{1,4}
\end{eqnarray*}
Using the above we can arrive at an expression for $U^{sol}_4$ in terms
of its class representative i.e $U_1$ as :
\begin{eqnarray*}
U_{1}A^{pth}_{1,4}U^*_4&=&B^{pth}_{1,4}\\
\Longrightarrow U_4&=& (B^{pth}_{1,4})^{-1} U_1 A^{pth}_{1,4}
\end{eqnarray*}
In general $U^{sol}=diag(U^{sol}_1,\cdots,U^{sol}_d)$ with $U^{sol}_i$ given by:
\begin{eqnarray*}
	U^{sol}_i=(B^{pth}_{c(i)i})^{-1} U_{c(i)} A^{pth}_{c(i)i}
\end{eqnarray*}
The next part of the proof is to show that $U^{sol}$ is indeed
Unitary, this follows because the non-zero sub-matrices are multiples of unitary as per 3) in Lemma 1. Suppose $A^*_{ij}A_{ij}=rI$ $r>0$, $A_{ij}$ can be written as
$\sqrt{r}$ x a Unitary matrix. By Lemma 1, 3) $B^*_{ij}B_{ij}$ is also $rI$, $B_{ij}$ is
also $\sqrt r$ x a Unitary matrix. We denote the product of all these
scalars ($\sqrt r$) along a path in $A^{pth}_{c(i)i}$ by $r^{pth}_{c(i)i}$.  Thus $A^{pth}_{c(i)i}$ is a multiple of
unitary given by $r^{pth}_{c(i)i}$ x a Unitary matrix, similarly $(B^{pth}_{c(i)i})^{-1}$ is
$(r^{pth}_{c(i)i})^{-1}$ x a Unitary matrix. Using these $U^{sol}_i$ (Definition 5) is: 
\begin{eqnarray*}
	U^{sol}_i&=&(r^{pth}_{c(i)i})^{-1} \times a \ Unitary \ matrix \times (r^{pth}_{c(i)i}) \times \ a \ Unitary \ matrix \\
	U^{sol}_i&=&(r^{pth}_{c(i)i})^{-1} \times (r^{pth}_{c(i)i}) \times \ a \ Unitary \ matrix \\
	U^{sol}_i&=& \ a \ Unitary \ matrix
\end{eqnarray*}
Hence $U^{sol}$ whose diagonal blocks are $U^{sol}_{i}$ is also Unitary.\\
The last part of the proof is to show that $U^{sol}A_lU^{sol*}=B_l$. This is same as showing $U^{sol}_iA_{ij}U^{sol*}=B_{ij}$ for $i,j \in \{1,\cdots,d\}$
and $l \in \{1,\cdots,p\}$.\\ 
Let $i,j$ belong to the same equivalence class in the $U$-induced partition of 
the given collection i.e $c(i)=c(j)$ by substituting for $U^{sol}_{i}$ and $U^{sol*}_j$ using Definition 5 and using Definition 4 for $pr(A_{ij})$ to simplify we have:
\begin{eqnarray*}
U^{sol}_iA_{ij}U^{sol*}_j&=&(B^{pth}_{c(i)i})^{-1}U_{c(i)}A^{pth}_{c(i)i} A_{ij}  (A^{pth}_{c(i)j})^{-1}U^*_{c(i)}B^{pth}_{c(i)j} \\
&=&(B^{pth}_{c(i)i})^{-1}U_{c(i)}[A^{pth}_{c(i)i} A_{ij} (A^{pth}_{c(i)j})^{-1}]U^*_{c(i)}B^{pth}_{c(i)j} \\
&=&(B^{pth}_{c(i)i})^{-1}U_{c(i)}pr(A_{ij})U^*_{c(i)}B^{pth}_{c(i)j} \\
\end{eqnarray*}
Further by Lemma 3 we have $U_{c(i)}pr(A_{ij})U^*_{c(i)}=pr(B_{ij})$, substituting in the above, using the the definition of $pr(B_{ij})$ we have :
\begin{eqnarray*}
U^{sol}_iA_{ij}U^{sol*}_j&=&(B^{pth}_{c(i)i})^{-1} pr(B_{ij}) B^{pth}_{c(i)j}   \\
&=&(B^{pth}_{c(i)i})^{-1} pr(B_{ij}) B^{pth}_{c(i)j}  \ \ \ (By \ Definition \ 4, of \ pr(B_{ij})) \\
&=& B_{ij}  \\
\end{eqnarray*}
In the above derivation $U_{c(i)}pr(A_{ij})U^*_{c(i)}=pr(B_{ij})$
is satisfied with $U_{c(i)}$ acting as a free variable ($\forall$ $U_{c(i)}$) since the given collection is in Solution form i.e 
$pr(A_{ij})$ , $pr(B_{ij})$ are multiples of identity. 

When $i,j$ do not beling to the same partition $A_{ij}=B_{ij}=0$. Hence the proof.
\subsection {Equivalent problem}
This section shows that when the collection is not in Solution-form we can set up an equivalent problem
such that we are a step closer to the simplest case. If the collection is not in Pre-Solution form this is done by diagonalizing submatrices along the diagonal i.e by diagonalizing $(A_{ii},B_{ii})$ or by diagonalizing $A_{ij}A_{ij}^*$
and $B_{ij}B_{ij}^*$ where $A_{ij},B_{ij}$ are off-diagonal sub-matrices.
On the other hand if the collection is in Pre-Solution form but not in Solution form it is done by diagonalizing $pr(A_{ij})$ and $pr(B_{ij})$ which are not multiples of identity. Here $pr(A_{ij})$ ($pr(B_{ij})$) is a multiple of Unitary and hence a Normal matrix.

We define an equivalent problem and relate it's solution to that of the original 
problem by Lemma 1.
 
\textbf{Definition 6:} Equivalent Problem.
Given a S.U.S problem with the collection $\{(\mathcal{A}_l,\mathcal{B}_l)\}_{l=1}^p$ where $\mathcal{A}_l,\mathcal{B}_l$ are $n$x$n$ Complex matrices and Unitary $Y,Z$, the S.U.S problem with $\{(Y\mathcal{A}_lY^*,Z\mathcal{B}_lZ^*)\}_{l=1}^p$ is an Equivalent problem of the given problem.

\textbf{Lemma 4:} Given an S.U.S problem and it's Equivalent problem
there exists a solution to the given problem i.f.f there is a solution to the original problem.

\textbf{Proof:} Let there exist a Unitary $U$ s.t $U\mathcal{A}_lU^*=\mathcal{B}_l$,
given a $U$ that solves the given problem, $\hat{U}=ZUY$ solves the Equivalent problem. Given $\hat{U}$ which solves the Equivalent problem $U=Z^*\hat{U}Y^*$ solves the given problem. Also when a sequence of Equivalent problems are constructed using $Y_1,\cdots,Y_t$ and $Z_1,\cdots,Z_t$ if the solution to the $t^{th}$ problem is known as $\hat{U}$ then the solution to the given problem can be obtained as $U=Z_1^* \cdots Z^*_t\hat{U}Y_t^* \cdots Y_1$.

The next Lemma helps to set up an equivalent problem based on the properties
of submatrices. 
 
\textbf{Lemma 5:} \begin{enumerate}
	\item A Complex square matrix $A$ is multiple of Identity if and only if both the Hermitian matrices $S_r=\frac{A+A^{*}}{2}$ and $S_c=\frac{A-A^{*}}{2j}$ are multiples of Identity.
	\item  If $A_{ij}$ is rectangular (non-square) $A_{ij}^*A_{ij}$ and $A_{ij}A^*_{ij}$ are both multiples of identity iff $A_{ij}=0$.
	\item When $A_{ij}$ is square, $A_{ij}^*A_{ij}$ and $A_{ij}A^*_{ij}$ are both multiples of identity iff $A_{ij}$ is multiple of Unitary.
\end{enumerate}

\textbf{Proof:} 
\begin{enumerate}
\item Suppose $A=\alpha I_n$, $\alpha=\alpha_r+j\alpha_c$ is a multiple of Identity then $S_r=\alpha_r I_n$ and $S_c=\alpha_c I_n$ are both multiples of Identity. Suppose $S_r=aI_n$ and $S_c=bI_n$ where $a,b$ are Real numbers, $A=S_r+jS_c$ hence $A=(a+jb)I$ which is multiple of Identity.
\item Given a non-square $A_{ij}$ of size $n_i$ x $n_j$, w.l.o.g let $n_i$ $>$  $n_j$ then the only way $A_{ij}A^*_{ij}$ can be a multiple of Identity is if $A_{ij}A^*_{ij}$ is zero as otherwise rank of $A_{ij}A^*_{ij}$ would be $n_i$. This is not possible since we know that rank of $A_{ij}A^*_{ij}$ is same as that of $A^*_{ij}A_{ij}$ which is $n_j$ and is $<$ $n_i$.
\item When $A_{ij}$ is square if $A^*_{ij}A_{ij}=rI$, $r>0$ then  
$\frac{A^*_{ij}}{\sqrt{r}}\frac{A_{ij}}{\sqrt r}=I$ which means $\frac{A_{ij}}{\sqrt r}$ is Unitary or $A_{ij}=\sqrt r$ x a Unitary matrix. If $r=0$ since rank of $A_{ij}$
is same as that of $A^*_{ij}A_{ij}$ it implies that $A_{ij}=0$x$I$. Hence the Proof.

\end{enumerate}
Theorem 2 is the main result of this section which sets up an Equivalent problem
when the collection is not in Solution form and ensures that the $U$ that solves the equivalent problem has more blocks than before thus moving a step closer to the
the simplest case which is that of a diagonal $U$.

\textbf{Theorem 2:} Given $U=diag(U_1,\cdots,U_d)$. If the given collection, ${\{(\mathcal{A}_l,\mathcal{B}_l)\}}_{l=1}^{p}$ is not in Solution form then
there is an Equivalent S.U.S problem with the collection ${\{(\hat{\mathcal{A}}_l,\hat{\mathcal{B}}_l)\}}_{l=1}^{p}$ s.t the unitary
$U$ that solves the Equivalent problem has the structure $U=diag(U_1,\cdots,U_d^{'})$ where $d'>d$.

\textbf{Proof:} Given the $U$-induced partition if the collection is not in Solution form it implies one of the following:
\begin{enumerate}
	\item There exists a $A_{ii}$ along the diagonal s.t $A_{ii}$ is not multiple of Identity.
	\item There exists a rectangular (non-square) $A_{ij}$ that is non-zero.
	\item There exists a off-diagonal square $A_{ij}$ which is not multiple of 
	Unitary.
	\item There exist $i,j$ belonging to class $c(i)$ s.t $pr(A_{ij})$ is not a
	multiple of Identity.
\end{enumerate}
In each of the above cases we show that we can find a Hermitian or a Normal matrix, diagonalize it and use the Unitary matrix that diagonalizes it to set up an Equivalent problem to solve in which $U$ has more diagonal blocks than $d$. 
Suppose 1. is true then by Lemma 5, 1. , one of the Hermitian matrices $\frac{A_{ii}+A_{ii}^*}{2}$ or $\frac{A_{ii}-A_{ii}^*}{2j}$ is not a multiple of Identity. Let us call this $n_i$x$n_i$ sized matrix $S$ ($n_i >1$). Diagonalizing $S$ we can write it as $D=Y_iSY_i^*$ where $Y_i$ is the Unitary matrix that diagonalizes $S$. \\
W.l.o.g let us assume $S=\frac{A_{ii}+A_{ii}^*}{2}$. If a solution were to exist to the S.U.S problem it is necessary that on $B$'s side $R=\frac{B_{ii}+B_{ii}^*}{2}$ be equal to $U_iSU_i^*$ and hence be unitarily similary to $S$. Let $Z_i$ diagonalize $R$ then $D=Z_iRZ_i^*$. \\
We use $Y=diag(I_1,\cdots,Y_i,\cdots,I_d)$ and $Z=diag(I_1,\cdots,Z_i,\cdots,I_d)$
to set up an equivalent problem as in Lemma 4. \\
Now to show that we have simplified the problem and the $U$ that solves the Equivalent problem has more blocks than $d$ let us look at the structure of $D$. Since $S,R$ are not multiples of identity they have atleast two distinct eigen values.
If the Equivalent problem were to have a solution then it is necessary that
$U_iDU_i^*=D$, assuming the eigen-values in $D$ are in descending order and using the
fact that any matrix commuting with a diagonal matrix has to be diagonal we have:
\begin{eqnarray*}	
	&U_iDU_i^*=D& , \ D=diag(\lambda_1I_{n_{\lambda_1}},\cdots,\lambda_{t}I_{n_{\lambda_t}}), \ t >1 \\
	&\Longrightarrow& U_iD=DU_i.\\
	&\Longrightarrow& U_i=diag(U^i_1,\cdots,U^i_t) \\
\end{eqnarray*}
Since this $U_i$ is one of the $d$ blocks of $U$, it now has $d'>d$
blocks since $t>1$.\\
In case we do not have Solution form because of 2. or 3. , by Lemma 5 either $A_{ij}A_{ij}^*$ or $A_{ji}^*A_{ji}$ will not be a multiple of Identity. Suppose $A_{ij}A_{ij}^*$ is not a multiple of Identity then using $Y_i$ and $Z_i$ that diagonalize $A_{ij}A_{ij}^*$ and  $B_{ij}B_{ij}^*$ respectively, as in the previous step we can diagonalize these with $Y_i,Z_i$ respectively and set-up an equivalent problem (Definition 6, Lemma 4) with $Y=diag(I_1,\cdots,Y_i,\cdots,I_d)$ and $Z=diag(I_1,\cdots,Z_i,\cdots,I_d)$.
On the other hand if $A_{ji}^*A_{ji}$ and $B_{ji}^*B_{ji}$ are not multiples of identity we can set up an equivalent problem using $Y_j$ and $Z_j$ that diagonalize 
$A_{ji}^*A_{ji}$ and $B_{ji}^*B_{ji}$ respectively with $Y=diag(I_1,\cdots,Y_j,\cdots,I_d)$ and $Z=diag(I_1,\cdots,Z_j,\cdots,I_d)$.
Rest of the argument to show $d'>d$ is as shown above for 1. \\
Suppose we do not have Solution form due to 4. let the class representative of $i,j$ be $c(i)$, by Lemma 3 w.k.t it is necessary that $U_{c(i)}pr(A_{ij})U^*_{c(i)}$ = $pr(B_{ij})$. Here
$pr(A_{ij})$ ($pr(B_{ij})$) is product of matrices which are multiples of Unitary and hence it is also a
multiple of Unitary. Since a matrix which is a multiple of Unitary is a Normal matrix ($mUm^*U^*=|m|^2I=m^*U^*mU$) it is unitarily diagonalizable. Just as in the case of 1. , 2. and 3. we can set up an equivalent problem using $Y_{c(i)}$ and $Z_{c(i)}$ which diagonalize $pr(A_{ij})$ and $pr(B_{ij})$ respectively.  Since $pr(A_{ij})$ and $pr(B_{ij})$ are not multiples of Identity we have $d'>d$.  Hence the proof.

Lemma 6 shows that the given collection is in Solution form when $U=diag(u_1,\cdots,u_n)$ is a diagonal matrix. This ensures that the algorithm runs for a maximum of $n$ steps.In the $n^{th}$ step the $U$ we are looking for will be diagonal and we end up in the simplest case.

\textbf{Lemma 6:} If the Unitary $U$ being solved for is a Diagonal matrix i.e $U=diag(u_1,\cdots,u_n)$, where $u_{i}$, $i \in \{1,\cdots,n\} $ are complex scalars with $\|u_i\|=1$, then the given collection is in Solution-form.

\textbf{Proof:}
 In $U$ all blocks are of size $1$ hence the submatrices in the $U$-induced partition of the collection are all $1$x$1$ matrices. The first of the conditions that diagonal blocks should be multiples of identity is trivially
satisfied as any $a_{ii}$ ($b_{ii}$) can be considered to be $a_{ii}$x$1$. There are no rectangular blocks hence the second condition is met. The off-diagonal $a_{ij}$'s ($b_{ij}$'s) are multiples of unitary since $a_{ij}=a_{ij}$x$1$ where $1$ is unitary. Hence the collection is in Pre-Solution
form. If we form the $U$-Induced Graph Partition and get $pr(A_{ij})$ and $pr(B_{ij})$ as per Definition 4 they will also be scalars and hence trivially multiples of identity.

Eventually after getting to Solution form if the given collection is S.U.S then $U^{sol}$ and all the $U$'s, $V$'s of the intermediate steps that diagonalize to set up equivalent problems can be combined as shown in Lemma 4 to get the solution to the original problem.

\subsection{Algorithm for S.U.Eq problem : Bi-Graph's Paths Partition $diag(U_1,\cdots,U_d)$ $\cup$ $diag(V_1,\cdots,V_f)$ }
The S.U.Eq problem which is to find if there exist Unitary $U,V$ s.t
s.t $U\mathcal{A}_lV^*=\mathcal{B}_l$ for a given collection of $m$x$n$ Complex matrices $\{(\mathcal{A}_l,\mathcal{B}_l)\}_{l=1}^{p}$ can also be solved by identifying the Solution form and the graph with paths that relates the blocks in $U=diag(U_1,\cdots,U_d)$ and in $V=diag(V_1,\cdots,V_f)$ to construct $U^{sol}$ and $V^{sol}$. Since the flow of the algorithm is similar to the description in \ref{subsec:flow} we describe the algorithm briefly highlighting the parts which are different.

\begin{figure}
	\centering
	\includegraphics[width=120mm]{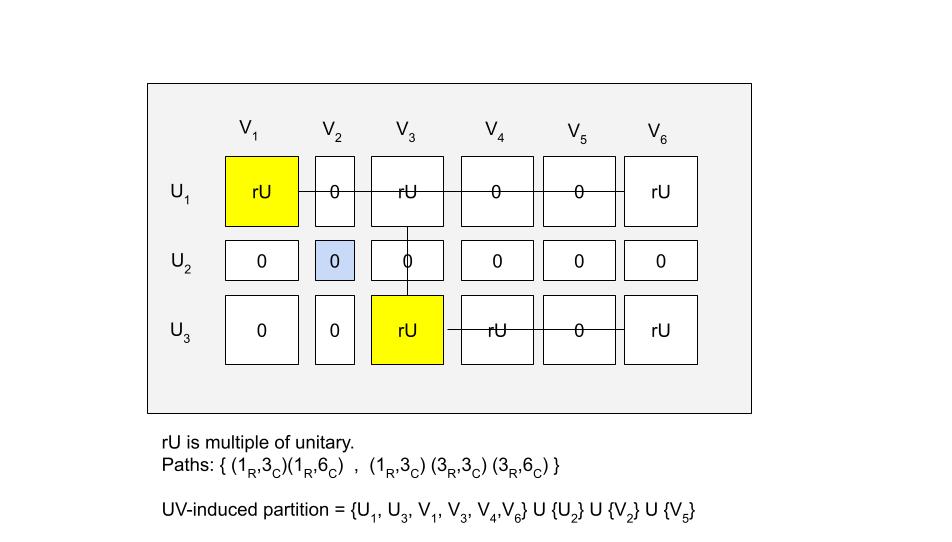}
	\caption{Pictorial depiction of S.U.Eq's Solution-form with paths and partition of blocks in $U$ and $V$ }
	\label{fig:fig2}
\end{figure}

In each iteration we solve for $U=diag(U_1,\cdots,U_d)$ and $V=diag(V_1,\cdots,V_f)$ with block sizes being $m_1,\cdots,m_d$ and $n_1,\cdots,n_f$ respectively.
The matrices of the collection are partitioned along the rows and columns based on these sizes to get the $U,V$-Induced Partition. Given this partition,
the collection is in Pre-Solution if the following is true:  
\begin{enumerate}
	\item Non-square submatrices are zeros i.e $A_{ij}=0=B_{ij}$.
	\item All off-diagonal square matrices are multiples of Unitary, i.e $A_{ij}A^*_{ij}=a_{ij}I=B_{ij}B^*_{ij}$, $a_{ij}\ge 0$
\end{enumerate}
The diagonal equations i.e $U_iA_{ii}V^*_i=B_{ii}$ involve both $U_i$ and $V_i$,
hence we do not have the `multiple of identity along the diagonal case' as in S.U.S. The graph associated with this Pre-Solution form is a bipartite graph (bi-graph) since the vertices of the row i.e $R=\{1_R,\cdots,d_R\}$ and the vertices of the columns i.e $C=\{1_C,\cdots,f_C\}$ are different. An $i\in R$ is connected by an edge to $j \in C$ if $A_{ij}$ ($B_{ij}$) is a non-zero multiple of Unitary. The paths in such a $U,V$-induced graph give the partition of $R \cup C$ (Figure \ref{fig:fig2}). Once we have the partition just as in the $S.U.S$ case, products can be associated with paths giving $A^{pth}_{c(i),i}$, $B^{pth}_{c(i),i}$, $pr(A_{ij})$, and $pr(B_{ij})$, $i,j \in R \cup C$. In Solution form with $pr(A_{ij})$, and $pr(B_{ij})$ being multiples of Identity, $U^{sol}_i$, $V^{sol}_j$ can be defined to get $U^{sol}$ and $V^{sol}$. \\
If the collection is not in Solution form an Equivalent problem is set up by diagonalizing either
$A_{ij}A^*_{ij}$ or $A^*_{ij}A_{ij}$ or $pr(A_{ij})$.  Suppose $A_{ij}A^*_{ij}$ is not multiple of Identity, let $Y_i$ be the matrix that diagonalizes
$A_{ij}A^*_{ij}$ and let $Z_i$ diagonalize $B_{ij}B^*_{ij}$ with the new collection
being $\hat{\mathcal{A}}_l$=$Y \mathcal{A}_l I$, $\hat{\mathcal{B}}_l$=$Z\mathcal{B}_l I$ where $Y=diag(I_1,\cdots,Y_i,\cdots,I_d)$ and $Z=diag(I_1,\cdots,Z_i,\cdots,I_f)$. If $A^*_{ij}A_{ij}$ is not a multiple of Identity then $\hat{\mathcal{A}}_l$=$ I \mathcal{A}_l Y$, $\hat{\mathcal{B}}_l$=$ I \mathcal{B}_l Z$ where  $Y=diag(I_1,\cdots,Y_j,\cdots,I_d)$ ,  $Z=diag(I_1,\cdots,Z_j,\cdots,I_f)$, $Y_j$ and $Z_j$ diagonalize $A^*_{ij}A_{ij}$
and $B^*_{ij}B_{ij}$ respectively. We can similarly get a simplified problem when $pr(A_{ij})$ is not a multiple of Identity. All these cases lead to an equivalent problem with $d'>d$ or $f'>f$. The problem gets solved in a maximum of $m+n$ steps.

\section{Computational Complexity}
\label{sec:comp_compl}
The flow of the algorithms as shown in \ref{subsec:flow} gives an idea of the computation involved in each iteration. The steps which may occur $n$ times in the worst case for the S.U.S problem are :
\begin{enumerate}
	\item Checking for Solution form.
	\item Setting up an Equivalent problem.
\end{enumerate}
When the collection is in Pre-Solution (Definition 2) form but not in Solution form (Definition 5) that is when
the maximum computation occurs. Assuming this happens every iteration gives us the worst-case complexity. The Complexity of partitioning the matrix along rows and
columns is $O(n^2)$ , for $p$ matrices it is $O(pn^2)$. Among the operations involving sub-matrices checking if a sub-matrix is multiple of Unitary is costliest and its cost can be estimated as follows:\\
There could be a maximum of $\frac{n}{n_i}$ x $\frac{n}{n_i}$ blocks of size $n_i$ for which we have to compute the product $A_{ij}A_{ij}^*$ whose cost is $O(n^3_i)$ ( $Kn^3_i$ ).  
The total cost would be:
\begin{eqnarray*}
	\Sigma_i \frac{n^2}{n^2_i}K n^3_i &=& K n^2 \Sigma_i n_i \\
	&=&K n^2 n \\
	&=&O(n^3)
\end{eqnarray*}
This repeated over maximum possible $n$ iterations for $p$ matrix tuples costs $O(pn^4)$.
The cost of forming the $U$-Induced Graph (Definition 3) is $O(pn^2)$. Once the graph is formed the cost of partitioning the graph by constructing the paths using either depth first or breadth first search is $O(pn^2)$ \citep{Trudeau1994grphtheory}. Over $n$ iterations this amounts to $O(pn^3)$. It is the cost of finding $A^{pth}_{c(i)i}$ and $pr(A_{ij})$ (Definition 4) which determines the worst case complexity. When finding $A^{pth}_{c(i)i}$ we can do so
incrementally so that for every new path the previously computed value is used and
only one additional product is computed, done this way the complexity of finding
$A^{pth}_{c(i)i}$ is as follows:\\
There are at most $\frac{n}{n_i}$ ($\Sigma_i n_i=n$) vertices corresponding to blocks of size $n_i$ each contributing a cost of $O(n^3_i)$. The total cost of finding $A^{pth}_{c(i)i}$ for all the vertices is:
\begin{eqnarray*}
	\Sigma_i \frac{n}{n_i}K n^3_i &=& K n \Sigma_i n^2_i \\
	&\le&K n  (\Sigma_i n_i)^2 \\
	&\le&K n n^2  \\
	&=&O(n^3)
\end{eqnarray*}
Over $n$ iterations this costs $O(n^4)$.
Finally $pr(A_{ij})$ ($pr(B_{ij})$ ) (Definition 4) is calculated using the $A^{pth}_{c(i)i}$ and
the complexity of this product is governed by number of possible blocks of size $n_i$ in the collection which is at most $\frac{pn^2}{n^2_i}$ and the 
cost of multiplication. These could be accounted for as follows:\\
\begin{eqnarray*}
	\Sigma_i \frac{pn^2}{n^2_i}K n^3_i &=& K pn^2 \Sigma_i n_i \\
	&=&K pn^2 \ n  \\
	&=&O(pn^3)
\end{eqnarray*}
Over $n$ iterations is $O(pn^4)$.
Thus the computation to check for Solution form
in the worst case has a complexity of $O(pn^4)$.\\
The next part is that of setting up of the equivalent problem which involves diagonalizing a Hermitian matrix or a Normal matrix (multiple of Unitary) followed by $p$ multiplications of the form $UA_iU^*$ ($VB_iV^*$). Diagonalization in the worst case costs $O(n^3)$ and $p$ multiplications cost $O(pn^3)$ hence the total cost of this step is $O(pn^3)$. Over $n$ steps it is $O(pn^4)$. \\
Finally once Solution form is reached
the cost of constructing $U^{sol}$ and finding $U$ is as follows:\\
$U^{sol}_{i}$ is calculated per vertex of the $U$-Induced Graph and hence there can be a maximum of $\frac{n}{n_i}$ of these of size $n_i$. Cost of $U^{sol}$ is:
\begin{eqnarray*}
	\Sigma_i \frac{n}{n_i}K n^3_i &=& K n \Sigma_i n^2_i \\
	&\le&K n n^2  \\
	&=&O(n^3)
\end{eqnarray*}
The cost of finding the Solution to the original problem using $U_{sol}$ is the computation $U=(\mathcal{V}_n\cdots \mathcal{V}_1)^*U_{sol}(\mathcal{U}_n\cdots \mathcal{U}_1)$ where $\{\mathcal{U}_i,\mathcal{V}_i\}_{i=1}^n$ are the unitary matrices that setup the equivalent problems. The total cost of this over $n$ iterations is of $O(n^4)$. Thus the cost of the algorithm we have presented to find a solution to the S.U.S problem is polynomial in $n,p$ and is of the order $O(pn^4)$.

Since the flow of the algorithm to solve the S.U.Eq problem is similar to that of the S.U.S solution similar arguments of worst case complexity can be made considering a maximum of $\frac{mn}{n^2_i}$ square blocks of size $n_i$ ($m_i$). The complexity of finding if the collection is in the Solution form which determines the overall complexity can be computed as follows:
\begin{eqnarray*}
	\Sigma_j \frac{mn}{n^2_j} K n^3_j &=& Kmn \Sigma_j  n_j  \\
	&=&Kmn^2 \\
	&=& O(mn^2)
\end{eqnarray*}
This repeated over a maximum of $m+n$ iterations for $p$ pairs is $O(p(m+n)mn^2)$ which is same as $O(pn^4)$ if we assume $m\approx n$.
   
\section{Implementation, Numerical Stability \& Accuracy, Canonical Features}
\label{subsec:numerical}
\subsection{Implementation}
Our python implementation of this algorithm to solve the S.U.S and S.U.Eq problems is available at https://github.com/harikrishna-vj/Simultaneous-Unitary-Similarity.git under the GNU Lesser General Public License. This implements the algorithm to solve the S.U.S problem as described in this paper. Using our implementation one could construct their own input collections and
check the outputs in each iteration. The $U$-Induced Partition, the paths,
the sub-matrices used to set up the Equivalent problem are all displayed 
after each iteration. Eventually if the given collection is S.U.S, $U^{sol}$ and the the $U$ that solves the problem is output. If the collection is not $S.U.S$ the reason for that, which is the mis-match in the eigen-values of the identified sub-matrices (Theorem 2) is output.
When the collection is S.U.S it is always interesting to check-out whether
the algorithm found a different $U$ than the one that was used to construct 
the problem ! 

It is challenging to construct examples that run for many iterations. Currently
three input collections are provided as examples to run the implementation. 
\begin{enumerate}
	\item A `$7$x$7$ example' with $p=3$, $2$-tuples which runs for $5$ iterations before hitting Solution-form and returing a $U$ that solves the
	problem.
	\item A `$5$x$5$ example' with $p=1$ in which $A,B$ are not Normal matrices.
	This is chosen since finding if a pair of non-normal matrices are unitarily similar is non-trivial. In $2$ iterations a $U$ that solves the problem is found. 
	\item Another `$5$x$5$ example' with $p=2$ in which the matrices of the collection are pair-wise unitarily similar is explored. The algorithm
	finds that the the collection is 'not S.U.S' in two iterations.
\end{enumerate}
We welcome those who are interested in this problem to explore the implementation, test it further by adding many more examples. 
The implementation can be be extended to add the algorithm for the S.U.Eq case,
implementing the bi-graph partition. Also the part to collect the 	`Canonical-features under Simultaneous Unitary Similarity' can be added to make it a complete implementation of the work described in this paper.

\subsection{Numerical Stability and Accuracy}
From the point of view of numerical stability and accuracy, the important step in each iteration of the algorithm is diagonalization of a Hermitian (Symmetric) matrix in case of the Pre-Solution form or diagonalization Normal matrix (multiple of Unitary) in the case of Solution form to get an equivalent problem. All the other computation is either matrix multiplication or checking for unitarity and equality. The diagonalization of Hermitian and Normal matrices is known to be numerically well-conditioned and stable. The error in eigen-value calculation ($\lambda_{err}$) of $S$ (matrix being diagonalized) given a machine precision of $\epsilon$ can be approximated to be:
\begin{equation}
 \lambda_{err} \le \epsilon || S ||
\end{equation}
The eigen vectors are sensitive to the `eigen-gap' i.e the smallest difference ($\delta$) between all eigen-value pairs of $S$. For Hermitian diagonalization, if $\theta$ is the angle between the actual and the estimated eigen-vectors then $sin( \theta$) as defined in Davis and Kahan's work \citep{Davis1970Sintheta}, for small $\theta$ can be approximated to be:
\begin{equation}
\sin(\theta) \le \frac{\epsilon|| S ||}{\delta}
\end{equation}
In our implementation we do the diagonalization at machine precision and comparisions required to establish orthogonality or equality at a precision which is few bits ($\Delta$) less than that of $\epsilon$. We sacrifice few bits to cover for the possible eigen-value error and eigen-vector sensitivity to eigen-gap. While this strategy is not fool-proof it gives a way to handle the accuracy issues which are inherent to the numerical problem of diagonalization.In our current implementation the comparisions are done at an accuracy of $10^{-9}$.
Improvements to the above implementation can look at tighter computation of error bounds \citep{Rump2023Fasterror}  adding pre-conditioning and dynamic multi-precision strategies as in \citep{Nicholas2025MixedPrecision}

\subsection{Classification under Simultaneous Unitary Similarity, A Canonical form}
Although our algorithms that solve the S.U.S and the S.U.Eq problems are for
a collection of $2$-tuples $\{(A_l,B_l)\}_{l=1}^p$ the steps can be executed for a
a collection of $p$ matrices $\{A_l\}_{l=1}^p$ noting important results in each step ending in Solution form and used as features of the
equivalance class to which $\{A_l\}_{l=1}^p$ belongs under Simultaneous Unitary Similarity.
Here we list the features that form part of the Canonical form for a given collection in each step. In every step the matrix partition is noted as 
a Canonical feature in the form $[n_1,n_2,\cdots,n_d]$. In the first iteration
it is $[n_1]$ where $n_1=n$.
If the collection is not in Solution form then as listed in the proof of Theorem 2
one of the following is true:
\begin{enumerate}
	\item There exists a $A_{ii}$ along the diagonal s.t $A_{ii}$ is not a multiple of Identity.\\
	Canonical feature: The matrix number $l$, $(i,i)$ indicating the submatrix and
	eigen-values of the Hermitian matrix $S$ (either $\frac{A+A^*}{2}$ or $\frac{A-A^*}{2j}$ ) that is diagonalized to setup the equivalent problem.
	
	\item There exists a rectangular (non-square) $A_{ij}$ that is non-zero.\\
	Canonical feature: The matrix number $l$ and $(i,j)$ indicating the submatrix and
	eigen-values of the Hermitian matrix $S$ (either $A^*_{ij}A_{ij}$ or $A_{ij}A^*_{ij}$) that is diagonalized to setup the equivalent problem.
	
	\item There exists a off-diagonal square $A_{ij}$ which is not multiple of Unitary.\\
	Canonical feature: The matrix number $l$, indices $i,j$ indicating the submatrix and
	eigen-values of the Hermitian matrix $S$ (either $A^*_{ij}A_{ij}$ or $A_{ij}A^*_{ij}$) that is diagonalized to set up the equivalent problem.
	
	\item There exist $i,j$ belonging to class $c(i)$ s.t $pr(A_{ij})$ is not a
	multiple of Identity.\\ 
	Here the collection is in Pre-Solution form hence the diagonal blocks are
	multiples of identity ($a_iI$) and the non-diagonals are multiples of unitary ($r_{ij}U$).\\
	Canonical feature: The $a_i$'s and the $r_{ij}$'s, the graph partition, The matrix number $l$, indices $i,j$ giving the $pr(A_{ij})$ which is not a multiple of identity,
	eigen-values of $S=pr(A_{ij})$ that is diagonalized to set up the equivalent problem.
\end{enumerate}
If the matrix is in Solution form then each $pr(A_{ij})$ is a multiple of identity $\beta_{ij}I$.
In Solution form the Canonical features are the $U$-Induced Graph, Partition and the $\beta_{ij}$'s. The procedure given here can be followed for the $S.U.Eq$ case as well to obtain the Canonical features under Simultaneous Unitary Equivalence. 
\section{Conclusion}
\label{sec:conclusion}
We have shown that the Simultaneous Unitary Similarity and the Simultaneous Unitary Equivalence problems can be solved by checking for a form which we refer to as the Solution form that captures the properties of the of matrices when we solve for a Unitary matrix that has the diagonal form. This is accomplished by constructing a solution i.e Unitary $U^{sol}$ ($V^{sol}$) using sub-matrices that are non-zero multiples of Unitary along the paths of a graph induced by the block diagonal structure of $U$. If in a iteration the collection is not in Solution we diagonalize a Hermitian matrix or a Normal matrix related to a submatrix of the partition and get a step close to the diagonal case. While the `eigen-gap' when diagonalizing is a factor in the numerical accuacy of the results the fact that we deal with mainly Hermitian matrices and Normal matrices for diagonalization ensures stability of the algorithm.  We also show that the steps of the algorithm if used on a given collection of $p$ matrices $\{\mathcal{A}_l\}_{l=1}^p$ until Solution form is reached retaining the Canonical features in each step gives a way to classify matrices upto Simultaneous Unitary Similarity. The time complexity of our algorithm is polynomial in $n,p$ and is of the order $O(pn^4)$. The S.U.S and S.U.Eq problems are well known for being computationally challenging.\citep{Segei1998wild}.

\bibliographystyle{unsrtnat}
\bibliography{references}  






\end{document}